\documentstyle[12pt]{article}
\begin{document}
\newtheorem{Def}{Definition}
\newtheorem{thm}{Theorem}
\newtheorem{lem}{Lemma}
\newtheorem{rem}{Remark}
\newtheorem{prop}{Proposition}
\newtheorem{cor}{Corollary}
\newtheorem{conj}{Conjecture}
\newtheorem{question}{Question}
\title
{Some Liouville theorems and applications}
\author{YanYan Li\thanks{Partially
 supported by
       NSF grant DMS-0401118.
To appear in a volume
in honor of Haim Brezis$'$  sixtieth birthday.}\\
\\
       Department of Mathematics\\
       Rutgers University\\
       110 Frelinghuysen Road\\
       Piscataway, NJ 08854\\
       USA
       \\
\\
\\
 {\it Dedicated to Haim
Brezis with  high respect  and friendship}
}
\date{}
\maketitle
\input { amssym.def}

\bigskip
\bigskip


\bigskip
\bigskip


\begin{abstract}  We give exposition of a Liouville theorem
established in \cite{Li3} which is a novel
 extension of the classical Liouville theorem
for harmonic functions. To illustrate
some ideas of the proof
of the  Liouville theorem, we present a new proof
 of the classical Liouville theorem
for harmonic functions.
Applications of the Liouville theorem, as well as that
of earlier ones in \cite{Li2},  can be found in
\cite{Li3, Li4} and \cite{W}.
\end{abstract}



The Laplacian operator $\Delta$ is invariant under
rigid motions:  For any function u on $\Bbb R^n$ and
for any rigid motion $T:
\Bbb R^n\to \Bbb R^n$,
$$
\Delta (u\circ T)= (\Delta u)\circ T.
$$

The following theorem  is classical:
\begin{equation}
u\in C^2, \ \
\Delta u=0\ \mbox{and}\ u>0\ \mbox{in}\ \Bbb R^n
\ \mbox{imply that}\ u\equiv \ \mbox{constant}.
\label{1}
\end{equation}

In this note we present a Liouville
theorem in \cite{Li3} which
is  a fully nonlinear version of
the classical Liouville theorem (\ref{1}).

Let  $u$  be a positive function in $\Bbb R^n$, and let
$\psi: \Bbb R^n\cup\{\infty\} \to \Bbb R^n \cup \{\infty\}$
be a M\"obius transformation, i.e. a transformation
generated by translations, multiplications by nonzero constants
and the inversion $x\to x/|x|^2$.  Set
$$
u_\psi:= |J_\psi|^{ \frac {n-2}{2n} } (u\circ \psi),
$$
where $J_\psi$ is the Jacobian of $\psi$.

It is proved in \cite{LL1} that an
operator $H(u, \nabla u, \nabla^2 u)$
is conformally invariant, i.e.
$$
H(u_\psi, \nabla u_\psi, \nabla ^2 u_\psi)
\equiv H(u, \nabla u, \nabla ^2u)\circ \psi
\ \mbox{holds for all positive $u$ and all M\"obius}\ \psi,
$$
 if and only if
$H$ is of the form
$$
H(u, \nabla u, \nabla^2 u)
\equiv f(\lambda(A^u))
$$
where
$$
A^u:= -\frac{2}{n-2}u^{  -\frac {n+2}{n-2} }
\nabla^2u+ \frac{2n}{(n-2)^2}u^ { -\frac {2n}{n-2} }
\nabla u\otimes\nabla u-\frac{2}{(n-2)^2} u^ { -\frac {2n}{n-2} }
|\nabla u|^2I,
$$
$I$ is the $n\times n$ identity matrix,
$\lambda(A^u)=(\lambda_1(A^u), \cdots, \lambda_n(A^u))$
 denotes the eigenvalues of $A^u$,
and $f$ is a function which is  symmetric  in $\lambda=(\lambda_1,
\cdots, \lambda_n)$.

Due to the above characterizing conformal
invariance property,   $A^u$ has been called
in the literature the conformal Hessian of $u$.
Since
$$
\sum_{i=1}^n \lambda_i(A^u)= -\frac 2{n-2} u^{ - \frac{n+2}{n-2} }
\Delta u,
$$
 Liouville theorem (\ref{1}) is equivalent to
\begin{equation}
u\in C^2, \ \ \lambda(A^u)\in \partial \Gamma_1
\ \mbox{and}\ u>0\ \mbox{in}\ \Bbb R^n
\ \mbox{imply that}\ u\equiv \ \mbox{constant},
\label{2}
\end{equation}
where
$$
\Gamma_1:= \{\lambda\ |\ \sum_{i=1}^n\lambda_i>0\}.
$$

Let
\begin{equation}
\Gamma\subset \Bbb R^n\
\mbox{be an open convex symmetric cone with vertex at the
origin}
\label{3}
\end{equation}
 satisfying
\begin{equation}
\Gamma_n:=\{\lambda\ |\ \lambda_i>0,
1\le i\le n\}\subset
\Gamma\subset \Gamma_1.
\label{4}
\end{equation}

Examples of such $\Gamma$ include those
given by elementary symmetric functions.  For
$1\le k\le n$, let
$$
\sigma_k(\lambda):=\sum_{1\le i_1< \cdots < i_k\le n}\lambda_{i_1}
\cdots \lambda_{i_k}
$$
be the $k-$th elementary symmetric function and let
$\Gamma_k:=\{\lambda\in \Bbb R^n\ |\
\sigma_1(\lambda), \cdots, \sigma_k(\lambda)>0\}$,
which is equal to the connected component of
$\{\lambda\in \Bbb R^n\ |\ \sigma_k(\lambda)>0\}$
containing the positive cone $\Gamma_n$,
 satisfies (\ref{3}) and (\ref{4}).

For an open subset $\Omega$ of $\Bbb R^n$, consider
\begin{equation}
\lambda(A^u)\in \partial \Gamma,\qquad \mbox{in}\ \Omega.
\label{5}
\end{equation}
The following definition of viscosity super and sub solutions of
(\ref{5}) has been given in \cite{Li3}.

\begin{Def}
A positive
continuous function $u$ in $\Omega$ is a viscosity subsolution
[resp. supersolution] of (\ref{5}) when the following holds: if
$x_0\in \Omega$, $\psi\in C^2(\Omega)$, $(u-\psi)(x_0)= 0$
and $u-\psi\le 0$ near
 $x_0$ then
$$
\lambda(A^\psi(x_0))\in \Bbb R^n\setminus \Gamma.
$$
[resp. if $(u-\psi)(x_0)=0$ and
 $u-\psi\ge 0$ near
 $x_0$   then
$\lambda(A^\psi(x_0))\in \overline \Gamma$].

We say that $u$ is a viscosity solution of (\ref{5}) if it is
both a viscosity supersolution and a viscosity  subsolution.
\end{Def}

\begin{rem}  If a positive $u$ in $C^{1,1}(\Omega)$ satisfies
$\lambda(A^u)\in \partial \Gamma$ a.e. in $\Omega$, then
it is a viscosity solution of (\ref{5}).
\end{rem}

Here is the Liouville theorem.

\begin{thm}\ (\cite{Li3})\ For $n\ge 3$, let
 $\Gamma$ satisfy (\ref{3}) and (\ref{4}), and let $u$ be a
positive locally Lipschitz viscosity solution of
\begin{equation}
\lambda(A^u)\in \partial \Gamma\qquad \mbox{in}\ \Bbb R^n.
\label{6} \end{equation}
Then $u\equiv u(0)$ in $\Bbb R^n$.
\label{thm-vis}
\end{thm}

\begin{rem} It  was proved
by  Chang, Gursky and Yang in
\cite{CGY} that
positive $C^{1,1}(\Bbb R^4)$ solutions to
$\lambda(A^u)\in \partial \Gamma_2$
are constants.
  Aobing Li proved in \cite{Lia} that
   positive $C^{1,1}(\Bbb R^3)$  solutions to
$\lambda(A^u)\in \partial \Gamma_2$
 are constants, and, for all $k$ and $n$,
   positive $C^3(\Bbb R^n)$
   solutions to
$\lambda(A^u)\in \partial \Gamma_k$  are constants.
The latter result for $C^3(\Bbb R^n)$
   solutions is independently established by
 Sheng, Trudinger and Wang
in \cite{STW}.
 Our proof is completely different.
\label{rem1}
\end{rem}

\begin{rem}
Writing $u=w^{ -\frac {n-2} 2}$, then
$$
A^u\equiv A_w:= w\nabla ^2 w -\frac 12 |\nabla w|^2 I.
$$
Theorem \ref{thm-vis}, with $\lambda(A^u)\in \partial\Gamma$
being replaced by $\lambda(A_w)\in \partial\Gamma$,
holds for $n=2$ as well.  See \cite{Li3}.
\end{rem}

In order to illustrate some of the ideas of our
proof of Theorem \ref{thm-vis} in \cite{Li3}, we give
a new proof of the classical Liouville theorem (\ref{1}).
We will derive (\ref{1}) by using the

\bigskip

\noindent {\bf  Comparison Principle for $\Delta$:}\
{\it Let $\Omega$ be a bounded open subset of $\Bbb R^n$ containing the
origin
$0$.  Assume that $u\in C^2_{loc}(\overline \Omega\setminus\{0\})$
and  $v\in C^2(\overline \Omega)$
satisfy
$$
\Delta u\le 0\ \ \mbox{in}\ \Omega\setminus\{0\}\qquad
\mbox{and}\qquad
\Delta v\ge 0\ \ \mbox{in}\ \Omega,
$$
and
$$
u>v\qquad\mbox{on}\ \partial \Omega.
$$
Then
$$
\inf_{  \Omega\setminus\{0\} }
(u-v)>0.
$$
}

\bigskip

It is easy to see from this proof of the Liouville theorem (\ref{1}) that
the
  following Comparison Principle for
locally Lipschitz  viscosity
solutions of (\ref{5}), established in \cite{Li2, Li3},
is sufficient for a proof of Theorem
\ref{thm-vis}.

\begin{prop}  Let $\Omega$ be a bounded open subset of $\Bbb R^n$
containing the origin $0$, and let
$u\in C^{0, 1}_{loc}(\overline \Omega\setminus\{0\})$
and $v\in C^{0, 1}(\overline \Omega)$.
Assume that
$u$  and $v$ are respectively positive  viscosity
supersolution and subsolution of (\ref{5}), and
$$
u>v>0\qquad \mbox{on}\ \partial \Omega.
$$
Then
$$
\inf_{\Omega\setminus\{0\}}(u-v)>0.
$$
\label{prop1}
\end{prop}

For the proof of Proposition \ref{prop1} and Theorem \ref{thm-vis},
see \cite{Li2, Li3}.
In this note, we give the

\medskip

\noindent{\bf Proof of
Liouville theorem  (\ref{1}) based on the
Comparison Principle for $\Delta$.}\
Let
$$
v(x):= \frac 12 [\min_{ |y|=1} u(y) ] |x|^{2-n},
\quad v_1(x):= \frac 1{|x|^{n-2}} v(\frac x{ |x|^2}),
\quad
u_1(x):= \frac 1{|x|^{n-2}} u(\frac x{ |x|^2}).
$$
Since $u_1$ and $v_1$ are still harmonic functions,
an application of the Comparison Principle for $\Delta$
on $\Omega:=$the unit ball yields
\begin{equation}
\liminf_{ |y|\to \infty} |y|^{n-2}u(y)>0.
\label{7}
\end{equation}

\begin{lem} For every $x\in \Bbb R^n$, there exists $\lambda_0(x)>0$ such that
$$
u_{x, \lambda}(y):= \frac{ \lambda^{n-2}}{  |y-x|^{n-2} }
u(x+\frac{  \lambda^2(y-x) }{  |y-x|^2  })\le u(y)\qquad
\forall\ 0<\lambda<\lambda_0(x), |y-x|\ge \lambda.
$$
\label{lem1}
\end{lem}

\noindent{\bf Proof.}\  Without loss of generality we may take $x=0$, and we
use  $u_\lambda$ to denote $u_{0,\lambda}$.
By the positivity and the Lipschitz regularity of $u$, there exists
$r_0>0$ such that
$$
r^{\frac{n-2}2} u(r, \theta)
< s^{\frac{n-2}2} u(s, \theta),\qquad\forall\
0<r<s<r_0, \ \theta\in \Bbb S^{n-1}.
$$
The above is equivalent to
\begin{equation}
u_\lambda(y)<u(y),\qquad 0<\lambda<|y|< r_0.
\label{8}
\end{equation}
We know from (\ref{7}) that, for some constant $c>0$,
$$
u(y)\ge c  |y|^{2-n},\qquad |y|\ge r_0.
$$
Let
$$
\lambda_0:= (\frac c{  \max_{ |z|\le r_0} u(z) })^{ \frac 1{n-2} }.
$$
Then
\begin{equation}
u_\lambda(y)\le (\frac{\lambda_0}{ |y|})^{n-2}
(\max_{ |z|\le r_0} u(z))\le  c |y|^{2-n} \le u(y),
\quad \forall\ 0<\lambda<\lambda_0, |y|\ge r_0.
\label{9}
\end{equation}
It follows from
(\ref{8}) and (\ref{9}) that
$$
u_\lambda(y)\le  u(y),
\quad \forall\ 0<\lambda<\lambda_0, |y|\ge \lambda.
$$
Lemma \ref{lem1} is established.

\vskip 5pt
\hfill $\Box$
\vskip 5pt

Because of Lemma \ref{lem1}, we may define, for any
$x\in \Bbb R^n$ and any $0<\delta<1$,
that
$$
\bar\lambda_\delta(x):=
\sup\{\mu>0\
|\ u_{x,\lambda}(y)\le (1+\delta)u(y),\
\forall\ 0<\lambda<\mu, |y-x|\ge \lambda\}\in (0, \infty].
$$

\begin{lem} For any
$x\in \Bbb R^n$ and any $0<\delta<1$,
$\bar \lambda_\delta(x)=\infty$.
\label{lem2}
\end{lem}

\noindent{\bf Proof.}\ We prove it by contradiction.  Suppose the
contrary, then, for some $x\in \Bbb R^n$ and some $0<\delta<1$,
$\bar \lambda_\delta(x)<\infty$. We  may assume, without loss
of generality,  that $x=0$, and  we use
$u_\lambda$ and $\bar \lambda_\delta$ to denote
respectively $u_{0, \lambda}$ and
$\bar \lambda_\delta(0)$.
Since the harmonicity is invariant under conformal
transformations and multiplication by constants, and since
$$
u(y)=
u_{\bar\lambda_\delta }(y)
< (1+\delta)u_{\bar\lambda_\delta}(y),
\qquad \forall\ |y|= \bar\lambda_\delta,
$$
an application of (\ref{7}) yields, using the fact
that $  (u_{\lambda})_\lambda\equiv u$,
$$
\inf_{ 0<|y|< \bar\lambda_\delta }
[(1+\delta)u_{\bar\lambda_\delta}(y) -u(y)]>0.
$$
Namely, for some constant $c>0$,
\begin{equation}
(1+\delta) u(y)- u_{\bar\lambda_\delta}(y)
\ge c|y|^{2-n},\qquad
\forall\ |y|\ge \bar\lambda_\delta.
\label{10}
\end{equation}
By the uniform continuity of
$u$ on the ball $\{ z\ |\ |z|<\bar \lambda_\delta\}$, there exists
$0<\epsilon<\bar \lambda_\delta$ such that
for all $ \bar \lambda_\delta\le \lambda\le \bar \lambda_\delta+\epsilon$,
and for all $|y|\ge \lambda$,
we have
\begin{eqnarray*}
(1+\delta)u(y)- u_\lambda(y)&\ge &
(1+\delta) u(y)- u_{\bar\lambda_\delta}(y)
+ [u_{\bar\lambda_\delta}(y) - u_\lambda(y)]
\\
&\ge & c|y|^{2-n} -
|y|^{2-n} |\lambda^{n-2} u(\frac{\lambda^2 y}{ |y|^2 })
- \bar \lambda_\delta^{n-2} u(\frac{ \bar \lambda_\delta^2 y}{ |y|^2 })|
\ge \frac c2|y|^{2-n}.
\end{eqnarray*}
This violates the definition of $\bar\lambda_\delta$.
Lemma \ref{lem2} is established.

\vskip 5pt
\hfill $\Box$
\vskip 5pt

By Lemma \ref{lem2},  $\bar\lambda_\delta\equiv \infty$
for all $0<\delta<1$.  Namely,
$$(1+\delta)u(y)\ge u_{x,\lambda}(y),
\quad \forall\ 0<\delta<1, x\in \Bbb R^n,
|y-x|\ge \lambda>0.
$$
Sending $\delta$ to $0$ in the above leads to
$$
u(y)\ge u_{x,\lambda}(y),
\quad \forall\  x\in \Bbb R^n,
|y-x|\ge \lambda>0.
$$
This easily implies $u\equiv u(0)$.
Liouville theorem (\ref{1}) is established.

\vskip 5pt
\hfill $\Box$

\end{document}